\documentclass[a4paper,12pt]{amsart}
\usepackage{amssymb}
\usepackage{ifthen}
\nonstopmode
\numberwithin{equation}{section}
\newtheorem{thm}{Theorem}[section]
\newtheorem{cor}[thm]{Corollary}
\newtheorem{lem}[thm]{Lemma}
\newtheorem{prop}[thm]{Proposition}

\theoremstyle{remark}
\newtheorem{rem}[thm]{Remark}
\newtheorem{example}[thm]{Example}

\newcounter{alphabet}
\newcounter{tmp}
\newenvironment{Thm}[1][]{\refstepcounter{alphabet}%
\bigskip%
\noindent%
{\bf Theorem \Alph{alphabet}}%
\ifthenelse{\equal{#1}{}}{}{ (#1)}%
{\bf .}
\itshape}{\vskip 8pt}

\newenvironment{pf}[1][]{%
 \vskip 3mm
 \noindent
 \ifthenelse{\equal{#1}{}}%
  {{\slshape Proof. }}%
  {{\slshape #1.} }%
 }%
{\qed\bigskip}

\newcounter{minutes}
\divide\time by 60
\newcounter{hours}
\multiply\time by 60
\addtocounter{minutes}{-\time}

\begin{document}
\newcommand{\A}{{\mathcal A}}
\newcommand{\B}{{\mathcal B}}
\newcommand{\T}{{\mathcal T}}
\newcommand{\M}{{\mathcal M}}
\newcommand{\F}{{\mathcal F}}
\newcommand{\Om}{{\Omega}}
\newcommand{\es}{{\mathcal S}}
\newcommand{\R}{{\mathbb R}}
\newcommand{\C}{{\mathbb C}}
\newcommand{\K}{{\mathcal K}}
\newcommand{\E}{{\mathcal E}}
\newcommand{\eK}{{\bold K}}
\newcommand{\uhp}{{\mathbb H}}
\newcommand{\Z}{{\mathbb Z}}
\newcommand{\N}{{\mathbb N}}
\newcommand{\D}{{\mathbb D}}
\newcommand{\UCV}{{\mathcal{UCV}}}
\newcommand{\kUCV}{{k\text{-}\mathcal{UCV}}}
\newcommand{\Perron}{{\mathcal P}}
\newcommand{\sphere}{{\widehat{\mathbb C}}}
\newcommand{\image}{{\operatorname{Im}\,}}
\renewcommand{\Im}{{\,\operatorname{Im}\,}}
\newcommand{\Aut}{{\operatorname{Aut}\,}}
\newcommand{\real}{{\,\operatorname{Re}\,}}
\renewcommand{\Re}{{\operatorname{Re}\,}}
\newcommand{\kernel}{{\operatorname{Ker}\,}}
\newcommand{\id}{{\operatorname{id}}}
\newcommand{\mob}{{\text{\rm M\"{o}b}}}
\newcommand{\Int}{{\operatorname{Int}\,}}
\newcommand{\Ext}{{\operatorname{Ext}\,}}
\renewcommand{\mod}{{\operatorname{mod}}}
\newcommand{\stab}{{\operatorname{Stab}}}
\newcommand{\SL}{{\operatorname{SL}}}
\newcommand{\PSL}{{\operatorname{PSL}}}
\newcommand{\PSU}{{\operatorname{PSU}}}
\newcommand{\tr}{{\operatorname{tr}}}
\newcommand{\diam}{{\operatorname{diam}\,}}
\newcommand{\inv}{^{-1}}
\newcommand{\area}{{\operatorname{Area}\,}}
\newcommand{\eit}{{e^{i\theta}}}
\newcommand{\eint}{{e^{in\theta}}}
\newcommand{\emint}{{e^{-in\theta}}}
\newcommand{\dist}{{\operatorname{dist}}}
\newcommand{\arctanh}{{\operatorname{arctanh}}}
\newcommand{\const}{{\operatorname{const.}}}
\newcommand{\capa}{{\operatorname{Cap}}}
\newcommand{\hdim}{{\operatorname{H-dim}}}
\newcommand{\rad}{{\operatorname{rad}}}
\newcommand{\partialb}{{\partial_{\operatorname{b}}}}
\newcommand{\CD}{{\operatorname{CD}}}
\newcommand{\hm}{{\mathcal H}}
\newcommand{\hc}{{\mathcal L}}
\newcommand{\cube}{{\mathcal Q}}
\newcommand{\Log}{{\,\operatorname{Log}}}

\newcommand{\Der}{{\mathfrak D}}
\newcommand{\X}{{\mathfrak X}}
\newcommand{\Isom}{{\operatorname{Isom}}}
\newcommand{\Hom}{{\operatorname{Hom}}}
\newcommand{\der}{{\mathcal D}}
\renewcommand{\r}{{\varphi}}
\newcommand{\s}{{\psi}}
\newcommand{\z}{{\partial/\partial z}}
\newcommand{\zf}{{\frac\partial{\partial z}}}
\newcommand{\zb}{{\partial/\partial \bar z}}
\newcommand{\zbf}{{\frac\partial{\partial \bar z}}}
\newcommand{\mone}{{\mbox{-}1}}
\newcommand{\iu}{{\operatorname{i}}}
\newcommand{\gs}{{\Sigma}}
\newcommand{\poly}{{\mathcal P}}
\newcommand{\V}{{\mathfrak D}}

\bibliographystyle{amsplain}
\title[Geometric invariants and univalence criteria]{
Geometric invariants associated with projective structures and
univalence criteria
}

\author[S.-A Kim]{Seong-A Kim}
\address{Department of Mathematics Education, Dongguk University \\
780-714, Korea}
\email{sakim@dongguk.ac.kr}  

\author[T. Sugawa]{Toshiyuki Sugawa}
\address{Graduate School of Information Sciences,
Tohoku University, Aoba-ku, Sendai 980-8579, Japan}
\email{sugawa@math.is.tohoku.ac.jp}

\begin{abstract}
For a nonconstant holomorphic map between projective Riemann surfaces
with conformal metrics, we consider invariant Schwarzian derivatives
and projective Schwarzian derivatives of general virtual order.
We show that these two quantities are related by
the ``Schwarzian derivative" of the metrics of the surfaces
(at least for the case of virtual orders $2$ and $3$).
As an application, we give univalence criteria for a meromorphic
function on the unit disk in terms of the projective Schwarzian
derivative of virtual order $3.$
\end{abstract}
\thanks{
The second author was supported in part by JSPS Grant-in-Aid for Scientific
Research (B), 17340039 and for Exploratory Research, 19654027.
}


\maketitle

\section{Introduction}

The (classical) Schwarzian derivative
\begin{equation}\label{eq:Sf}
Sf=\left(\frac{f''}{f'}\right)'-\frac12\left(\frac{f''}{f'}\right)^2
=\frac{f'''}{f'}-\frac32\left(\frac{f''}{f'}\right)^2
\end{equation}
of a nonconstant meromorphic function $f$ on a plane domain
was introduced by Schwarz to construct a conformal mapping
of a Jordan domain bounded by finitely many circular arcs.
The reason why the Schwarzian derivative is so useful is that
it satisfies the invariance relation
$S(M\circ f\circ L)=Sf\circ L\cdot(L')^2$
for M\"obius transformations $L$ and $M.$
In particular, the quantity $\gs f=\lambda^{-2}Sf$ for
a function on the unit disk $\D=\{z\in\C: |z|<1\}$ is invariant
in the sense that $\gs(M\circ f\circ T)=(\gs f)\circ T\cdot (T'/|T'|)^2$
for a M\"obius transformation $M$ and an analytic automorphism
$T$ of $\D.$
Here, $\lambda(z)|dz|=|dz|/(1-|z|^2)$ is the hyperbolic metric of $\D.$
Due to these invariance properties, the Schwarzian derivative
has found many applications
in complex analysis, Teichm\"uller theory, 1-dimensional dynamical systems,
and so on.

It is thus a natural desire to seek for more quantities analogous to the
Schwarzian derivative.
Indeed, Schwarzian derivatives of higher order were proposed
in \cite{Ahar69}, \cite{Tamanoi96} and \cite{Schip00}.
Those Schwarzians certainly enjoy several interesting properties
but they do not find many applications so far.
One reason perhaps comes from the lack of invariance.
For instance, as we will see in Section 5,
Schwarzians of Aharonov and Tamanoi are invariant
under post-composition with M\"obius transformations but not
under pre-composition (as differentials), in general.

The authors proposed in \cite{KS:isd} the
{\it invariant} Schwarzian derivative
$\gs^n f$ of virtual order $n$ for a nonconstant holomorphic map
$f$ from a Riemann surface with conformal metric into another.
This derivative satisfies the invariance relation 
$\gs^n(M\circ f\circ T)=(\gs^n f)\circ T\cdot (T'/|T'|)^n$
for local isometries $M$ and $T$ (see Lemma \ref{lem:inv} for a more
precise formulation) at the expense of analyticity.
It involves conformal metrics of both the source and target
surfaces and therefore has a complicated form in general.
Note that $\gs^2f$ is nothing but the above $\gs f$ when
the source and target surfaces are $\D$ and $\sphere=\C\cup\{\infty\}$
equipped with the hyperbolic and spherical metrics, respectively.

In the present paper, we will introduce yet another kind of
Schwarzian derivatives, denoted by $V^nf$ and called the {\it projective}
Schwarzian derivative of virtual order $n,$
for a nonconstant holomorphic map $f$ from a projective surface with
a conformal metric into a projective surface.
Here, we note that projective structure is finer than complex structure
and that a plane domain has a natural projective structure
(see the next section for details).
This satisfies the invariance property that $V^n(M\circ f\circ T)
=(V^nf)\circ T\cdot (T')^n$ for a projective map $M$ and
a projective local isometry $T$ (see Lemma \ref{lem:V} below for details).
Since $V^nf$ does not involve a conformal metric of the target surface,
the form of $V^nf$ is much simpler than that of $\gs^nf.$

One of the most important applications of the Schwarzian derivative
is a univalence criterion due to Nehari \cite{Nehari49}.

\begin{Thm}[Nehari]\label{Thm:Nehari}
Let $f$ be a nonconstant meromorphic function on the unit disk $\D.$
If $f$ is univalent in $\D,$ then $\|Sf\|_2\le 6.$
Conversely, if $\|Sf\|_2\le 2,$ then $f$ must be univalent in $\D.$
The numbers $6$ and $2$ are both sharp.
\end{Thm}

Here, we set
$$
\|\varphi\|_c=\sup_{z\in\D} (1-|z|^2)^{c}|\varphi(z)|
$$
for a $\sphere$-valued function $\varphi$ on $\D$ and a real number $c.$
We notice the invariance property that
\begin{equation}\label{eq:inv}
\|\varphi\circ T\cdot|T'|^c\|_c=\|\varphi\|_c
\end{equation}
holds for each analytic automorphism $T$ of $\D$ because
of the formula $|T'(z)|(1-|z|^2)=1-|T(z)|^2.$

The first assertion in Theorem A
was indeed found by Kraus as early as in 1932, and re-discovered by Nehari
later.
Therefore, it is sometimes called the Kraus-Nehari theorem.
Theorem A constitutes a basis of 
the theory of Teichm\"uller spaces.
See \cite{Lehto:univ} for details.

A similar result for the pre-Schwarzian derivative $f''/f'$ is also known
and it is utilized to construct another model of the universal Teichm\"uller
space (see \cite{AG86} or \cite{Sug07ut}).

In the present paper, as a by-product of our investigation,
we give a univalence criterion for a function $f$ on $\D$
in terms of the projective Schwarzian
\begin{equation}\label{eq:Vf}
Vf(z)=V^3f(z)=(Sf)'(z)-\frac{4\bar z}{1-|z|^2}Sf(z)
\end{equation}
of virtual order $3.$

\begin{thm}\label{thm:main}
Let $f$ be a nonconstant meromorphic function on the unit disk $\D.$
If $f$ is univalent in $\D,$ then $\|Vf\|_3\le 16.$
The number $16$ is sharp.
Conversely, if $\|Vf\|_3\le 3/2,$ then $f$ is univalent in $\D.$
\end{thm}

It appears that the constant $3/2$ in the theorem is not sharp.
On the other hand, the constant cannot be replaced by a number greater
than $16\sqrt3/9\approx3.0792$ as we will see in Example \ref{ex:S}.
In the proof, we will see that the theorem is not stronger than the Nehari
univalence criterion in Theorem A.
We, however, expect that this quantity $Vf$ would open a new window
to a family of univalence criteria, as a paper by Duren, Shapiro and Shields
\cite{DSS66} led to Becker's univalence criterion \cite{Becker72}.

Let us summarize the contents of the paper.
We recall the definition and basic properties of the Schwarzian derivative
of a holomorphic map between projective (Riemann) surfaces in Section 2.
In Section 3, we review basics of the Peschl-Minda derivatives
and Schwarzian derivatives of higher order due to Aharonov \cite{Ahar69},
Tamanoi \cite{Tamanoi96} and the authors \cite{KS:isd}
for a nonconstant holomorphic map 
between Riemann surfaces with conformal metrics.
Section 4 is devoted to a relation between $\gs f$ and $Sf$
when the surfaces are projective and have conformal metrics
(Theorem \ref{thm:isd}).
To this end, we introduce the Schwarzian derivative of a conformal metric.
This result has several applications as we will see there.

The higher-order Schwarzians of Aharonov and Tamanoi
cannot be extended to holomorphic maps between projective Riemann surfaces
unlike the classical Schwarzian.
In Section 5, we introduce {\it projective} Schwarzian derivatives.
We will then try to generalize Theorem \ref{thm:isd} 
for projective Schwarzians of order 3.
Our future task is to extend this result to the case of general order.
The last section will be devoted to the proof of Theorem \ref{thm:main}
and to the computation for a concrete example.

\medskip
\noindent
{\bf Acknowledgments.}
The second author presented a talk based on this research
at a seminar held in W\"urzburg on November 2008.
The authors are grateful for useful comments to the audience,
especially, Richard Fournier, Daniela Kraus, Oliver Roth,
Stephan Ruscheweyh, and Vagia Vlachou.

\section{Projective structures}

Let us briefly recall basic properties of the Schwarzian derivative
$Sf,$ given in \eqref{eq:Sf},
of a nonconstant meromorphic function $f.$
It is well known that $Sf\equiv 0$ if and only if $f$ is (a restriction of)
a M\"obius transformation and that the formula
\begin{equation}\label{eq:cayley}
S(g\circ f)=(Sg)\circ f\cdot(f')^2+Sf
\end{equation}
holds for the composite map $g\circ f.$
In particular, $S(M\circ f\circ L)=Sf\circ L\cdot(L')^2$
for M\"obius transformations $L$ and $M$ as we already mentioned in Introduction.

We are tempted to define the Schwarzian derivative for a holomorphic
map $f$ between Riemann surfaces.
The above formula, however, tells us that the value $Sf$ may depend on the
choice of local coordinates.
Hence, we are naturally led to the idea to restrict the local coordinates
so that the transition functions are M\"obius, that is, the idea
of projective structures.

Though the notion of projective structures is standard, we describe
its basics in some detail in order to clarify the formulations below.
The notion of projective structures is obtained by replacing holomorphic
maps by M\"obius maps in the definition of complex structure
(see \cite[\S 9]{Gunning:R}).
More precisely, a {\it projective structure} on a surface $R$ is the
equivalence class of an atlas
$\{z_\alpha:U_\alpha\to U_\alpha'\}_{\alpha\in A},$ where
$U_\alpha\subset R,~U_\alpha'\subset\C$ are open sets and
$z_\alpha:U_\alpha\to U_\alpha'$ is a homeomorphism for $\alpha\in A$
such that the transition function $z_\beta\circ z_\alpha^{-1}$ is a M\"obius
map on each connected component of $z_\alpha(U_\alpha\cap U_\beta)$
for $\alpha,\beta\in A.$
Such an atlas will be called a projective atlas.
Two atlases are defined to be equivalent if the union of the two
is again a projective atlas.

A {\it projective surface} will mean a surface with a projective structure.
A map $z_\alpha:U_\alpha\to U_\alpha'$ in a projective atlas of a
projective surface $R$ will be called a {\it projective coordinate.}
Note that a projective surface has the canonical complex structure,
which is called the underlying complex structure.
In other words, a projective structure is finer than a complex structure,
and thus a projective surface can be regarded as a Riemann surface
in a canonical way.
A continuous map $f$ from a projective surface $R$ into another
projective surface $R'$ is called
{\it projective} if $w\circ f\circ z\inv$ is either M\"obius or constant
whenever $z$ and $w$ are projective coordinates of $R$ and $R',$ respectively.

For instance, a plane domain $\Omega$ has the natural atlas
$\{\id:\Omega\to\Omega\},$ which gives rise to a projective structure
on $\Omega.$
In the sequel, unless otherwise stated, a plane domain will be endowed
with this {\it natural} projective structure.
The uniformization theorem states that the universal covering surface
of a Riemann surface $R$ is conformally equivalent to one (and only one)
of the standard
surfaces; the Riemann sphere $\sphere=\C\cup\{\infty\}=\C_{+1},$
the complex plane $\C=\C_{0},$
and the unit disk $\D=\{z\in\C: |z|<1\}=\C_{-1}.$
Here the notation $\C_\delta,~\delta=+1,0,-1,$
is introduced to handle with these three at once.
According to the cases $\delta=+1,0,-1,$ the surface $R$ is called
elliptic, parabolic, or hyperbolic, respectively.
Let $h:\C_\delta\to R$ be a holomorphic universal covering projection
of $\C_\delta$ onto $R.$
Since the group of conformal automorphisms of $\C_\delta$ is
contained in the group of M\"obius transformations,
the local inverses of $h$ give rise to a projective structure on $R.$
This projective structure will be called {\it standard.}
The standard projective structure is characterized by the property that
a holomorphic universal covering projection of the standard domain $\C_\delta$
onto $R$ is projective.

For projective coordinates $z$ of $R$ and $w$ of $R',$ we define
a meromorphic quadratic differential on $R$ by
$S_{R,R'}f=S(w\circ f\circ z\inv)dz^2$ for a nonconstant
holomorphic map $f:R\to R'.$
Then the meromorphic quadratic differential $S_{R,R'}f=S_{R,R'}(z)dz^2$
does not depend on the choice of the projective coordinates,
thus it is well defined.
In other words, the system of functions
$S_{R,R',\alpha}f=S(w\circ f\circ z_\alpha\inv)$ for projective coordinates
$z_\alpha$ of $R$ and $w$ of $R',$ we have the relations
$S_{R,R',\beta}f\circ g_{\beta,\alpha}\cdot(g_{\beta,\alpha}')^2
=S_{R,R',\alpha}f$ for $g_{\beta,\alpha}=z_\beta\inv\circ z_\alpha.$
Note that $S_{R,R'}f\equiv0$ if and only if $f$ is a nonconstant projective map.
If $R$ and $R'$ are plane domains (with natural projective structures),
then obviously $S_{R,R'}f$
coincides with the usual Schwarzian derivative $Sf(z)dz^2.$
If we do not need to refer to the projective structures of $R$ and $R',$
we write $S_{R,R'}f=Sf$ simply.
For basic information about projective structures and applications to
Teichm\"uller spaces, see \cite{Nag:teich} and references therein.

\section{Invariant Schwarzian derivative}

We first recall the Peschl-Minda derivatives.
See \cite{KS07diff} or \cite{Schip07} for details.

For the sake of simplicity, we bigin with the case of plane domains.
Let $\Omega$ and $\Omega'$ be plane domains with conformal metrics
$\rho=\rho(z)|dz|$ and $\sigma=\sigma(w)|dw|,$ respectively.
Throughtout the present paper, a conformal metric will be always smooth.

The $\rho$-derivative of a smooth function $\r$ on $\Omega$ is defined by
$$
\partial_\rho\r
=\frac{\partial\r}{\rho}
=\frac1{\rho(z)}\frac{\partial\r(z)}{\partial z}.
$$
For a holomorphic map $f:\Omega\to\Omega',$ we define
the Peschl-Minda derivative
$D^nf$ of order $n$ with respect to $\rho$ and $\sigma$
inductively by
\begin{align*}
D^1f&=\frac{\sigma\circ f}{\rho}f' \\
D^{n+1}f&
=\left[\partial_\rho-n\partial_\rho(\log\rho)
+(\partial_\sigma\log\sigma)\circ f\cdot 
D^1f\right]D^nf \quad(n\ge1).
\end{align*}
We further set
$$
Q^nf=\frac{D^{n+1}f}{D^1f},\quad n\ge1.
$$
Since $D^nf$ and $Q^nf$ depend on the metrics, 
we often write $D_{\sigma,\rho}^nf$ and $Q_{\sigma,\rho}^nf$ for them.
The following formula will be useful below.

\begin{lem}[\cite{KS:isd}]\label{lem:dQ}
$$
\partial_\rho(Q^nf)=Q^{n+1}f-\big[Q^1f-n\partial_\rho\log\rho\big] Q^nf.
$$
\end{lem}

We next recall the definitions of Aharonov invariants and Tamanoi's
Schwarzian derivatives.
Let $f$ be a meromorphic function on a domain $D$ in the complex plane.
For $z\in D$ with $f(z)\ne\infty, f'(z)\ne0,$ we expand
$$
\frac{f'(z)}{f(z+w)-f(z)}
=\frac1w-\sum_{n=1}^\infty \psi_n[f](z)w^{n-1}
$$
for small enough $w.$
The quantities $\psi_n[f](z)$ were introduced by Aharonov \cite{Ahar69}
and called the {\it Aharonov invariants} by Harmelin \cite{Har82}.
Independently, Tamanoi \cite{Tamanoi96} defined 
the Schwarzian derivative $S_n[f]$ of virtual order $n$ by
$$
\frac{f'(z)(f(z+w)-f(z))}{\frac12f''(z)(f(z+w)-f(z))+f'(z)^2}
=\sum_{n=0}^\infty S_n[f](z)\frac{w^{n+1}}{(n+1)!}.
$$
Note that $S_0[f]=1$ and $S_1[f]=0.$
The Aharonov invariants and Tamanoi's Schwarzian derivatives are essentially
same in the sense that the following relation holds:
$$
\sigma_n[f]=\psi_n[f]+\sum_{k=2}^{n-2}\psi_k[f]\sigma_{n-k}[f],
\quad n\ge 2,
$$
where $\sigma_n[f]=S_n[f]/(n+1)!.$
In particular, $3!\psi_2[f]=S_2[f]=Sf$ and $4!\psi_3[f]=S_3[f]=(Sf)'.$
We also have the following recursive relations:
\begin{equation}\label{eq:rec}
S_{n}[f]
=S_{n-1}[f]'+\tfrac12 S_2[f]\sum_{k=1}^{n-1} \binom{n}{k}S_{k-1}[f] S_{n-k-1}[f],
\quad n\ge 3
\end{equation}
(see \cite{KS:isd} for details).
For example, $S_4[f]=S_3[f]'+4S_2[f]^2.$

We note that $S_n[M\circ f]=S_n[f]$ holds for a M\"obius transformation $M.$
However, unlike $Sf=S_2[f],$ the higher-order Schwarzian derivative
$S_n[f]$ does not behave nicely with pre-composition with 
M\"obius transformations.
For instance, $S_3[f\circ L]=(S_2[f\circ L])'=S_3[f]\circ L\cdot (L')^3
+2S_2[f]\circ L\cdot L'L''$ for a M\"obius transformation $L.$
Therefore, in general, $\sigma_n[f]dz^n$ and $S_n[f]dz^n$ 
are not invariant under the change of projective coordinates.
Moreover, it is unlikely that a result similar to the second half of
Theorem A holds for these Schwarzians.
For instance, we have no constant $c\ge0$ such that $\|S_3[f]\|_3\le c$ 
implies univalence of $f.$
Indeed, the function $f(z)=e^{az}$ is not univalent in $\D$ for $|a|>\pi$
but $S_3[f]=0.$

In the sequel, we mainly consider Tamanoi's Schwarzian derivatives.
We now give another description of them.
Define a sequence of polynomials $P_n=P_n(x_1,\dots,x_n)$ of 
$n$ indeterminates
$x_1,\dots,x_n$ inductively by $P_0=1, P_1=0, P_2=x_2-3x_1^2/2,$ and
$$
P_{n}=\sum_{k=1}^{n-1}(x_{k+1}-x_1x_k)\frac{\partial P_{n-1}}{\partial x_k}
+\frac12 P_2\sum_{k=1}^{n-1} \binom{n}{k}P_{k-1}P_{n-k-1},\quad n\ge3.
$$
For instance, $P_3=x_3-4x_1x_2+3x_1^3.$
Then, by letting $q_n[f]=f^{(n+1)}/f',$ we have (see \cite{KS:isd})
$$
S_n[f]=P_n(q_1[f],q_2[f],\dots,q_n[f]),\quad n\ge0.
$$

By using the above expression of $S_n[f],$
we define the higher-order invariant Schwarzian derivatives
$\gs^nf$ for a nonconstant
holomorphic map $f:\Omega\to \Omega'$ between plane domains $\Omega$ and $\Omega'$
with conformal metrics $\rho$ and $\sigma,$ respectively, by
$$
\gs^nf=P_n(Q^1f,\dots,Q^nf),\quad n\ge0.
$$
We sometimes write $\gs_{\sigma,\rho}^nf$ for $\gs^nf$ to indicate the
metrics involved.
For brevity, we also write $\gs^2f=\gs f.$
More concretely,
\begin{equation}\label{eq:gs2}
\gs f
=Q^2f-\frac32 (Q^1f)^2
=\frac{D^3f}{D^1f} - \frac{3}{2}\left(\frac{D^2f}{D^1f}\right)^2.
\end{equation}

We can deduce the following formula from \eqref{eq:rec} (see \cite{KS:isd}):

\begin{equation}\label{eq:rec2}
\gs^{n}f=\big(\partial_\rho-(n-1)\partial_\rho\log\rho\big)\gs^{n-1}f
+\frac12\gs^2f\sum_{k=1}^{n-1} \binom{n}{k}\gs^{k-1}f \gs^{n-k-1}f
\end{equation}
for $n\ge3.$

We record an invariance property of these quantities in the following form.

\begin{lem}\label{lem:inv}
Let $\Omega, \hat\Omega, \Omega', \hat\Omega'$ 
be plane domains with conformal
metrics $\rho, \hat\rho, \sigma, \hat\sigma,$ respectively.
Suppose that locally isometric holomorphic maps
$g:\hat\Omega\to \Omega$ and $h:\Omega'\to\hat\Omega'$ are given.
Then, for a nonconstant holomorphic map $f:\Omega\to\Omega',$
the following transformation rule is valid:
\begin{align*}
Q^n_{\hat\sigma,\hat\rho}(h\circ f\circ g)
&=(Q^n_{\sigma,\rho}f)\circ g
\cdot\left(\frac{g'}{|g'|}\right)^n \\
\gs^n_{\hat\sigma,\hat\rho}(h\circ f\circ g)
&=(\gs^n_{\sigma,\rho}f)\circ g
\cdot\left(\frac{g'}{|g'|}\right)^n.
\end{align*}
\end{lem}

See \cite[Lemma 3.6]{KS07diff}) for the proof of the relation for $Q^n.$
For the proof of the relation for $\gs^n,$ we observe that
$P_n$ is of weight $n,$ that is,
$P_n$ is a linear combination of monomials of weight $n.$
Here, the weight of a monomial $x_{j_1}\cdots x_{j_k}$ is defined to be
the number $j_1+\cdots+j_k.$
Therefore, $\gs^nf$ obeys the same transformation rule as that of $Q^n f.$
See \cite[Lemma 4.1]{KS:isd} for details.

By the last lemma, $Q^nf$ and $\gs^nf$
can be defined for a nonconstant holomorphic map $f$ between
Riemann surfaces with conformal metrics as a suitable differential form.
In particular, if $f:R\to R'$ is a locally isometric holomorphic map,
then the local coordinates can be chosen so that $f=\id$ and $\sigma=\rho,$
and hence, $D^nf=0$ for $n\ge2$ and $ Q^nf=\gs^nf=0$ for $n\ge1.$
Note that the Riemann surfaces are {\it not} required to be projective
for the definition of $Q^nf$ and $\gs^nf.$

\section{Relationship between invariant and classical Schwarzians}

It is fundamental to have a relation between the invariant Schwarzian
$\gs f=\gs^2f$ and the classical Schwarzian $Sf.$
To this end, we introduce a few quantities associated with a conformal metric.

Let $\rho=\rho(z)|dz|$ be a (smooth) conformal metric on a Riemann surface $R.$
We recall that the Gaussian curvature of $\rho$ is defined as
\begin{equation}\label{eq:kappa}
\kappa_\rho=-\frac{\Delta\log\rho}{\rho^2}
=-4\frac{\partial\bar\partial\log\rho}{\rho^2}.
\end{equation}
Note that $\kappa_\rho$ does not depend on the particular choice of
local coordinates.

For the hyperbolic metric $\lambda_\Omega$ of a plane domain $\Omega,$
Minda \cite{Minda97} considered the Schwarzian derivative
$2\partial^2\log\lambda_\Omega-2(\partial\log\lambda_\Omega)^2.$
We can define the same quantity for a conformal metric on a projective surface.
\begin{lem}\label{lem:Theta}
For a projective surface $R$ with a conformal metric $\rho,$ let
\begin{align*}
\Theta_{R,\rho}(z)
&=
2\frac{\partial^2\rho}{\rho}-4\left(\frac{\partial\rho}{\rho}\right)^2
=2\partial^2\log\rho-2(\partial\log\rho)^2 \\
&=2\left(\zf\right)^2\log\rho(z)
-2\left(\frac{\partial\log\rho(z)}{\partial z}\right)^2.
\end{align*}
Then $\Theta_{R,\rho}(z)dz^2,$
evaluated with projective coordinates $z,$ 
are pieced together to a smooth quadratic
differential $\Theta_{R,\rho}$ on $R.$
\end{lem}

\begin{rem}
Here and in the sequel, we adopt this redundant-looking notation
$\Theta_{R,\rho}$ because this quantity depends not only on the metric $\rho,$
but also on the projective structure of $R.$
Thus the reader should note that $R$ in the subscript notation indicates
rather the projective structure than the underlying surface.
\end{rem}

\begin{pf}
Let $\hat z$ be another projective coordinate of $R$ and let
$g=\hat z\circ z^{-1}$ be the transition function, namely,
$\hat z=g(z).$
If we write $\rho=\rho(z)|dz|=\hat\rho(\hat z)|d\hat z|,$
then  $\rho(z)=\hat\rho(g(z))|g'(z)|.$
Thus, as in \cite{Minda97}, 
\begin{equation}\label{eq:rho}
\partial\log\rho=(\partial\log\hat\rho)\circ g\cdot g'
+\frac{g''}{2g'}
\end{equation}
and therefore,
\begin{equation*}
2\partial^2\log\rho
-2\left(\partial\log\rho\right)^2
=2\left[(\partial^2\log\hat\rho)\circ g-(\partial\log\hat\rho)^2\circ g
\right](g')^2+Sg.
\end{equation*}
Since $g$ is M\"obius, $Sg=0.$ The proof is now complete.
\end{pf}

The following well-known fact is important in the sense that
the metric of constant curvature yields a holomorphic quadratic differential.
By solving the Schwarzian differential equation, in turn, we can reproduce
the metric (see \cite{KR:cm} for details).
For convenience of the reader, we supply a short proof as well.

\begin{lem}\label{lem:c}
Let $R$ be a projective surface with smooth conformal metric $\rho.$
Then $\Theta_{R,\rho}$ is a holomorphic
quadratic differential on $R$ if and only if $\rho$ has constant
Gaussian curvature.
\end{lem}

\begin{pf}
By \eqref{eq:kappa}, we have $\bar\partial\partial\log\rho=-\kappa_\rho\rho^2/4$
and therefore 
$\bar\partial\Theta_{R,\rho}
=-\rho^2\partial\kappa_\rho/2.$
Thus we see now that $\Theta_{R,\rho}$ is holomorphic if and only if
$\partial\kappa_\rho=0.$
Since $\kappa_\rho$ is real-valued, the last condition is equivalent
to that $\kappa_\rho$ be constant.
\end{pf}

\begin{example}\label{ex:standard}
The standard domain $\C_\delta$ has the complete
metric $\lambda_\delta=(1+\delta|z|^2)\inv|dz|$
of constant Gaussian curvature $4\delta$ 
and has the natural projective structure.
These metrics are called {\it spherical, Euclidean,} and {\it hyperbolic}
according to the cases when $\delta=+1,0,$ and $-1.$
The metric $\lambda_\delta$ will be called {\it standard.}
Since $\partial\log\lambda_\delta=-\delta\bar z/(1+\delta|z|^2),$
we easily see that $\Theta_{\C_\delta,\lambda_\delta}=0$
for $\delta=+1,0,-1.$
\end{example}

Let $h:\C_\delta\to R$ be a holomorphic universal covering projection
of $\C_\delta$ onto a Riemann surface $R.$
Since the covering transformations are contained 
in the group $\Isom^+(\C_\delta)$ of isometries on $(\C_\delta,\lambda_\delta),$
the metric $\lambda_\delta$ projects to a metric
$\lambda_R,$ which will be called the standard metric of $R.$
Thus, $\lambda_R$ is a smooth conformal metric on $R$ of constant Gaussian
curvature $4\delta$ such that $h^*(\lambda_R)=\lambda_\delta.$
Here we record the following observation.

\begin{lem}\label{lem:standard}
Let $R$ be a Riemann surface with standard metric $\lambda_R$ and
standard projective structure.
Then $\Theta_{R,\lambda_R}=0.$
\end{lem}

The following result connects the invariant Schwarzian derivative $\gs f$ 
with the classical one $Sf.$

\begin{thm}\label{thm:isd}
Let $R,R'$ be projective surfaces with smooth conformal metrics
$\rho,\sigma,$ respectively, and let $f:R\to R'$ be a nonconstant
holomorphic map.
Then
\begin{equation}\label{eq:gs}
\gs_{\rho,\sigma} f
=\rho^{-2}\left[S_{R,R'}f +f^*\Theta_{R',\sigma}-\Theta_{R,\rho}\right],
\end{equation}
where $f^*\Theta_{R',\sigma}$ is the pull-back 
$(\Theta_{R',\sigma}\circ f)(f')^2$
as a quadratic differential.
\end{thm}

\begin{pf}
By taking projective coordinates, we may assume that $R$ and $R'$ are
plane domains.
First, by Lemma \ref{lem:dQ}, we have
$$
Q^{2}f=\partial_\rho(Q^1f)+\big[Q^1f-\partial_\rho\log\rho\big] Q^1f.
$$
Substituting the last formula to \eqref{eq:gs2}, we also have
$$
\gs f=\partial_\rho Q^1f-\frac12(Q^1f)^2-(\partial_\rho\log\rho)Q^1f.
$$
Since
\begin{align}\notag
Q^1f&=\frac{D^2f}{D^1f}=2\partial_\rho\log(\sigma\circ f)+\partial_\rho\log f'
-2\partial_\rho\log\rho \\
&=\rho\inv[2(\partial\log\sigma)\circ f\cdot f'+f''/f'-2\partial\log\rho],
\label{eq:Q1}
\end{align}
the relation
\begin{align*}
\partial_\rho Q^1f
&=
-(\partial_\rho\log\rho)Q^1f \\
&\quad +\rho^{-2}[2(\partial^2\log\sigma)\circ f\cdot (f')^2
+2(\partial\log\sigma)\circ f\cdot f''
+(f''/f')'-2\partial^2\log\rho]
\end{align*}
holds.
Therefore, we have
\begin{align*}
\rho^2\gs f&=
2(\partial^2\log\sigma)\circ f\cdot (f')^2
+2(\partial\log\sigma)\circ f\cdot f''+(f''/f')'-2\partial^2\log\rho \\
&\quad-\frac12(\rho Q^1f)^2-2(\partial\log\rho)(\rho Q^1f).
\end{align*}
In view of \eqref{eq:Q1} we compute
\begin{align*}
&~~\frac12(\rho Q^1f)^2+2(\partial\log\rho)(\rho Q^1f) 
=\frac12(\rho Q^1f+4\partial\log\rho)(\rho Q^1f) \\
&=\frac12\big[2(\partial\log\sigma)\circ f\cdot f'+f''/f'\big]^2
-2(\partial\log\rho)^2 \\
&=2(\partial\log\sigma)^2\circ f\cdot(f')^2
+2(\partial\log\sigma)\circ f\cdot f''+\frac12(f''/f')^2
-2(\partial\log\rho)^2,
\end{align*}
and substitute it to the last expression of $\rho^2\gs f$
to get the required relation.
\end{pf}

\begin{cor}
Let $R$ be a Riemann surface with conformal metric $\rho$
and let $\varepsilon\in\{+1,0,-1\}.$
For a holomorphic map $f:R\to\C_\varepsilon$ and a M\"obius
transformation $M$ with $M(f(R))\subset\C_\varepsilon,$
$$
\gs_{\rho,\lambda_\varepsilon}^n(M\circ f)
=\gs_{\rho,\lambda_\varepsilon}^n f,
\quad n\ge0.
$$
\end{cor}

\begin{pf}
In view of the formula \eqref{eq:rec2}, it suffices to show
the relation for $n=0,1,2$ by induction.
The relation trivially holds for $n=0,1.$
Thus we may assume that $n=2.$
We assign a projective structure (e.g., the standard one) to $R$ so that
we regard $R$ as a projective surface.
Since $\Theta_{\C_\varepsilon,\lambda_\varepsilon}=0$
(see Example \ref{ex:standard}), we have
$$
\gs^2(M\circ f)
=\rho^{-2}[S_{R,\C_\varepsilon}(M\circ f)-\Theta_{R,\rho}]
=\rho^{-2}[S_{R,\C_\varepsilon}f-\Theta_{R,\rho}]
=\gs^2 f
$$
by Theorem \ref{thm:isd}.
\end{pf}

Also, by Example \ref{ex:standard}, we obtain the following.

\begin{cor}
Let $\delta,\varepsilon\in\{-1,0,+1\}.$
For a nonconstant holomorphic map $f:\C_\delta\to\C_\varepsilon,$ the
following relation holds:
$$
\gs f=\lambda_\delta^{-2}Sf.
$$
\end{cor}

This relation was observed from time to time for various combinations
of standard metrics, see \cite{KM93}, \cite{KM01} and \cite{MM97two}.
Note that a nonconstant holomorphic map $f:\C_\delta\to\C_\varepsilon$
exists if and only if $\delta\le\varepsilon.$

Let $R$ be a projective surface with standard metric $\lambda_R.$
The quadratic differential $\theta_R=\Theta_{R,\lambda_R}$ is called
the {\it uniformizing connection} of a projective surface $R$
(see \cite{Kra89} for the case when $R$ is a hyperbolic plane domain).
We remark that, when the projective structure of $R$ is standard,
we have $\theta_R=0$ by Lemma \ref{lem:standard}.
The following result gives a way of computing the universal covering
projection of a given surface $R$ once we have an explicit form of
$\theta_R.$
This idea traces back to Poincar\'e.

\begin{cor}
Let $R$ be a projective surface and let $h:\C_\delta\to R$ be a
holomorphic universal covering projection.
Then the uniformizing connection $\theta_R$ is a holomorphic
quadratic differential on $R$ and related to $h$ by
$$
h^*\theta_R\equiv(\theta_R\circ h)\cdot(h')^2=-Sh.
$$
\end{cor}

\begin{pf}
By Lemma \ref{lem:c}, we can see that $\theta_R$ is holomoprhic.
Note that $\Theta_{\C_\delta,\lambda_\delta}=0$ by Example \ref{ex:standard}
and that $\gs h=0$ because $h$ is a local isometry.
We now apply Theorem \ref{thm:isd} to obtain 
$S_{\C_\delta,R}h+h^*\theta_R=0,$
which is nothing but the required relation.
\end{pf}

\section{Projective Schwarzian derivatives of higher order}

As we noted, Tamanoi's Schwarzian derivatives are not well defined
for holomorphic maps between projective Riemann surfaces.
However, if the source surface is equipped with a conformal metric,
it is possible to define another sort of Schwarzian derivatives of higher order.
We begin by recalling a differential-geometric tool to do so.
Let $\varphi=\varphi(z)dz^n$ be a smooth $n$-differential on a Riemann
surface $R$ with conformal metric $\rho.$
Then
$$
\Lambda_\rho(\varphi)
=\big[\partial\varphi-2n(\partial\log\rho)\varphi\big]dz^{n+1}
$$
is a well-defined $(n+1)$-differential on $R.$
Indeed, $\Lambda_\rho(\varphi)$ is nothing but the covariant derivative
of $\varphi dz^n$ in $z$-direction
with respect to the Levi-Civita connection of $\rho$
(see \cite[\S 3]{KS07diff} for details).

Based on \eqref{eq:rec}, we can express $S_n[f]$
in terms of $Sf$ and its higher derivatives in the same way as before.
Define a sequence of polynomials $T_n=T_n(x_2,\dots, x_n)$ of $n-1$ 
indeterminates with integer coefficients, inductively, by $T_2=x_2$ and
$$
T_n=\sum_{k=2}^{n-1}\frac{\partial T_{n-1}}{\partial x_k}\cdot x_{k+1}
+\frac{x_2}2\sum_{k=1}^{n-1}\binom{n}{k}T_{k-1}T_{n-k-1},\quad n\ge 3.
$$
Here, we also set $T_0=1$ and $T_1=0.$
For instance, $T_3=x_3,~T_4=x_4+4x_2^2$ and $T_5=x_5+13x_2x_3.$
We can also easily verify that $T_n$ is of weight $n.$
Then, we have
\begin{equation}\label{eq:der}
S_n[f]=T_n(Sf,(Sf)',\dots, (Sf)^{(n-2)}),\quad n\ge 3.
\end{equation}

Let $R$ and $R'$ be projective surfaces and let $\rho$ be a conformal
metric on $R.$
For a nonconstant holomorphic map $f:R\to R',$
we define differentials $\V_{R,\rho,R'}^nfdz^n~(n\ge 2)$ on $R$
inductively by $\V_{R,\rho,R'}^2f=S_{R,R'}f$ and
$$
\V_{R,\rho,R'}^nfdz^n=\Lambda_\rho(\V_{R,R'}^{n-1}fdz^{n-1})\quad(n\ge 3),
$$
where differentiations are performed with respect to projective coordinates.
Namely, $\V_{R,\rho,R'}^nfdz^n=\Lambda_\rho^{n-2}(S_{R,R'}fdz^2).$
Furthermore, we define $V_{R,\rho,R'}^nf$ by
$$
V_{R,\rho,R'}^nf=T_n(\V_{R,\rho,R'}^2f,\dots,\V_{R,\rho,R'}^nf)
$$
for $n\ge 2.$
Here, the product in the substitution of $\V_{R,\rho,R'}^k$'s
is understood as the tensor product.
Since $T_n$ is of weight $n,$ $V_{R,\rho,R'}^nf$
can be regarded as an $n$-differential on $R.$
We call $V_{R,\rho,R'}^nf$ the {\it projective Schwarzian derivative}
of virtual order $n$ for a nonconstant holomorphic map $f:R\to R'.$
When we do not need to indicate the projective structures and/or
the conformal metric,
we simply write $\V^nf$ or $\V_{\rho}^nf$ for $\V_{R,\rho,R'}^nf.$
We do the same thing for $V.$
Note that $V^2f$ is nothing but the classical Schwarzian $Sf,$
which is independent of the metric $\rho.$
Furthermore, we have
$$
V^3f
=S_3[f]-4\frac{\partial\rho}\rho S_2[f]
=(Sf)'-4\frac{\partial\rho}\rho Sf
$$
and
$$
V^4f=S_4[f]-10\frac{\partial\rho}\rho S_3[f]+4\left[
7\left(\frac{\partial\rho}\rho\right)^2-\frac{\partial^2\rho}\rho
\right]S_2[f].
$$
If $R$ and $R'$ are plane domains with the Euclidean metric $|dz|,$
we obviously have $\V^nf=(Sf)^{(n-2)}$ and, by \eqref{eq:der},
$V^nf=S_n[f]$ for $n\ge2.$

For $R=\C_{-1}$ and $R'=\C_{+1}$
with standard metrics and standard projective structures,
$\V^3f=V^3f$ is same as in \eqref{eq:Vf}.

It is convenient for future reference to rephrase explicitly the fact that
$V^nf$ is a {\it well-defined} $n$-differential
in the following way.

\begin{lem}\label{lem:V}
Let $\Omega$ be a plane domain with conformal metric $\rho.$
For a M\"obius transformation $g,$ set
$\hat\Omega=g\inv(\Omega)$ and $\hat\rho=g^*\rho.$
For a nonconstant meromorphic function $f$ on $\Omega$
and a M\"obius transformation $h,$ the following relation holds:
$$
V_{\hat\rho}^n(h\circ f\circ g)=(V_{\rho}^nf)\circ g\cdot(g')^n,
\quad n\ge 2.
$$
\end{lem}

Since every analytic automorphism of $\C_{-1}=\D$ is M\"obius and isometric
with respect to $\lambda_{-1},$ the following result can be derived immediately
(see also \eqref{eq:inv}).

\begin{cor}\label{cor:inv}
Let $f$ be a nonconstant meromorphic map on $\D.$
For an analytic automorphism $T$ of $\D$ and a M\"obius
transformation $M,$
$$
V^n(M\circ f\circ T)=V^nf\circ T\cdot(T')^n.
$$
In particular, $\|V^n(M\circ f\circ T)\|_n=\|V^nf\|_n,~ n\ge 2.$
Here $V^n$ is defined for the hyperbolic metric on $\D.$
\end{cor}

We also need to consider the ``derivatives" of $\Theta_{R,\rho}.$
For a projective surface $R$ with conformal metric $\rho=\rho(z)|dz|,$ set
$\Theta_{R,\rho}^ndz^n=\Lambda_\rho^{n-2}(\Theta_{R,\rho}dz^2)$
for $n\ge2.$
It is a basic problem to find a relation between $\gs^nf$ and $V^nf$
for a nonconstant holomorphic map $f$ between projective surfaces
with conformal metrics.
We treat, however, with the case when $n=3$ only.
Compare with Theorem \ref{thm:isd}.

\begin{thm}\label{thm:V}
Let $R,R'$ be projective surfaces with smooth conformal metrics
$\rho,\sigma,$ respectively, and let $f:R\to R'$ be a nonconstant
holomorphic map.
Then
$$ 
\gs_{\sigma,\rho}^3f
=\rho^{-3}\left[V_{R,\rho,R'}^3f +f^*\Theta^3_{R',\sigma}
-\Theta_{R,\rho}^3\right]
+2\rho^{-2}f^*\Theta_{R',\sigma}^2Q_{\sigma,\rho}^1f,
$$
where $f^*\Theta_{R',\sigma}^n$ is the pull-back 
$(\Theta_{R',\sigma}^n\circ f)(f')^n$
as an $n$-differential.
\end{thm}

\begin{pf}
By \eqref{eq:rec} with $n=3,$ we have
$$
\gs^3f=\partial_\rho\gs f-2\partial_\rho\log\rho\cdot\gs f.
$$
Letting $U=\rho^2\gs f,$ we now see that
$$
\rho^3\gs^3f
=\rho^3\partial_\rho(\rho^{-2}U)-2\rho\partial_\rho\log\rho\cdot U
=\partial U-4\partial\log\rho\cdot U.
$$
Substitution of \eqref{eq:gs} into the last formula yields
\begin{align*}
\rho^3\gs^3f&=
[(Sf)'+(\partial\Theta_{R',\sigma})\circ f\cdot(f')^3
+\Theta_{R',\sigma}\circ f\cdot(2f'f'')-\partial\Theta_{R,\rho}] \\
&\quad
-4\partial\log\rho[Sf +\Theta_{R',\sigma}\circ f\cdot(f')^2-\Theta_{R,\rho}] \\
&=V_{R,\rho,R'}^3 f+f^*\Theta_{R',\sigma}^3-\Theta_{R,\rho}^3 \\
&\quad +\left[4(\partial\log\sigma)\circ f\cdot(f')^3+2f'f''
-4(\partial\log\rho)(f')^2\right]\Theta_{R',\sigma}\circ f.
\end{align*}
The required relation now follows from \eqref{eq:Q1}.
\end{pf}

In conjunction with Lemma \ref{lem:standard}, we have the following.

\begin{cor}
For a nonconstant holomorphic map $f:\C_\delta\to\C_\varepsilon,$ the
following relation holds:
$$
\gs^3f=\lambda_\delta^{-3}V_{\lambda_\delta}^3f.
$$
\end{cor}

\section{Applications to univalence criteria}

Aharonov \cite{Ahar69} (see also \cite{Har82})
showed that a nonconstant meromorphic function
$f$ on the unit disk $\D$ is univalent if and only if
$$
\sum_{n=1}^\infty n\left|
\sum_{k=1}^n \binom{n-1}{k-1}(-\bar z)^{n-k}(1-|z|^2)^{k+1}
\psi_{k+1}[f](z)\right|^2\le1,\quad z\in\D.
$$
In particular, the first term gives the Kraus-Nehari theorem
(the first half of Theorem A.
The second term gives the inequality
$$
\big|(1-|z|^2)^3\psi_3[f](z)-\bar z(1-|z|^2)^2\psi_2[f](z)\big|
=\frac{(1-|z|^2)^3|Vf(z)|}{24}\le \frac1{\sqrt2}
$$
for $f$ univalent in $\D.$
Thus we have $\|Vf\|_3\le 12\sqrt2\approx 16.97,$
which is slightly worse than the estimate in Theorem \ref{thm:main}.
Here, $Vf=V_{\lambda_{-1}}^3f$ with the notation introduced in the
previous section.

We are now in a position to prove Theorem \ref{thm:main}.
First let $f$ be univalent in $\D.$
For an arbitrary $a\in\D,$ letting $T(z)=(z+a)/(1+\bar az),$
by Corollary \ref{cor:inv} we have the relation
$$
(1-|a|^2)^3|Vf(a)|=|(Vf)\circ T(0)||T'(0)|^3=|V(f\circ T)(0)|.
$$
Since $\|V(f\circ T)\|_3=\|Vf\|_3$ by Corollary \ref{cor:inv},
it is enough to show that $|Vf(0)|\le 16$
for univalent meromorphic function $f$ on $\D.$
We may further assume that $f(0)\ne\infty.$
Then we look at the Laurent expansion
$$
\frac{f'(0)}{f(w)-f(0)}
=\frac1w-\sum_{n=0}^\infty c_nw^{n}
=\frac1w-\sum_{n=1}^\infty \psi_n[f](0)w^{n-1}
$$
to have $\psi_3[f](0)=c_2.$
It is known that $|c_2|\le 2/3$ and equality holds if and only if
$f'(0)/(f(w)-f(0))=(1-e^{i\theta}w^3)^{2/3}/w$ for a real constant $\theta$
(see, for instance, \cite[Theorem 4.6, p.~135]{Duren:univ}).
Therefore, we have
$$
|Vf(0)|=|S_3[f](0)|=24|\psi_3[f](0)|=24|c_2|\le 16.
$$
The sharpness 
is also clear from the above argument.

The latter part of Theorem \ref{thm:main} follows from the next proposition
together with Theorem A.

\begin{prop}\label{prop:main}
For a locally univalent meromorphic function $f$ on the unit disk,
the inequalities
$$
\frac{16}{25\sqrt5}\|Vf\|_3
\le\|Sf\|_2\le \frac43\|Vf\|_3
$$
hold. Here, the constant $16/25\sqrt5$ is sharp.
\end{prop}

We are left to prove Proposition \ref{prop:main} only.
The following representation formula will be the main tool for the proof.

\begin{lem}\label{lem:int}
Let $f$ be a locally univalent meromorphic function on the unit disk
with $\|Vf\|_3<\infty.$
Then
$$
Sf(z)=-\frac1{\pi}\iint_{|\zeta|<1}
\frac{(1-|\zeta|^2)^4Vf(\zeta)}{(1-|z|^2)^4(\bar\zeta-\bar z)}
d\xi d\eta \quad(\zeta=\xi+i\eta).
$$
\end{lem}

\begin{pf}
Let $\rho(z)=1/(1-|z|^2)$ as before.
Then, by definition, we have
\begin{equation}\label{eq:de}
Vf=(Sf)'-4\frac{\partial\rho}{\rho}Sf
=\rho^4\partial(\rho^{-4}Sf)
\quad\text{on} \D.
\end{equation}
Define functions $\psi$ and $\varphi$ on $\C$ by 
$\psi=\rho^{-4}Sf,~\varphi=\rho^{-4}Vf$ on $\D$
and $\psi=\varphi=0$ on $\C\setminus\D.$
Since $\rho^{-3}Vf$ is bounded on $\D$ by assumption,
we see that $\varphi\in C_0^0(\C),$
where $C_0^k(\C)$ denotes the class of continuous functions on $\C$ 
with compact support which have continuous partial derivatives up to
order $k.$

We show now that $\psi\in C_0^1(\C).$
Consider the function $y(t)=Sf(t\zeta)$ in $0\le t<1$ for a fixed 
$\zeta\in\partial\D.$
By \eqref{eq:Vf}, we see that $y(t)$ satisfies the linear differential equation
$$
y'(t)-\frac{4t}{1-t^2}y(t)=v(t),
$$
where $v(t)=\zeta\cdot Vf(t\zeta).$
Thus
$$
(1-t^2)^2y(t)=\int_0^t (1-s^2)^2v(s)ds+y(0).
$$
Since $|v(s)|\le \|Vf\|_3(1-s^2)^{-3},$ we obtain
$$
(1-t^2)^2|y(t)|\le 
\|Vf\|_3\int_0^t(1-s^2)^{-1}ds+|Sf(0)|
=\frac12\|Vf\|_3\log\frac{1+t}{1-t}+|Sf(0)|.
$$
In view of \eqref{eq:Vf}, we thus have
$$
Sf(z)=O\left(\frac1{(1-|z|^2)^2}\log\frac{1+|z|}{1-|z|}\right)
\quad\text{and}\quad
(Sf)'(z)=O\left(\frac1{(1-|z|^2)^3}\log\frac{1+|z|}{1-|z|}\right)
$$
as $|z|\to 1^-.$
In particular, $\rho(z)^{-3}Sf(z)\to0$ and $\rho(z)^{-4}(Sf)'(z)\to0$
when $|z|\to1^-.$
It is now easy to verify that $\psi\in\C_0^1(\C).$
Hence, by \eqref{eq:de} it is confirmed that $\partial\psi=\varphi$ on $\C$
and, equivalently, $\bar\partial\bar\psi=\bar\varphi$ on $\C.$

Consider now the integral transform
$$
h(z)=-\frac1\pi\iint_\C\frac{\overline{\varphi(\zeta)}}{\zeta-z}d\xi d\eta
$$
of $\bar\varphi.$
It is known (see \cite{Ahlfors:qc}) that $h$ is H\"older continuous
with exponent $<1$ and satisfies 
$\bar\partial h=\bar\varphi$ (in the sense of distribution).
Therefore, $\bar\partial(h-\bar\psi)=\bar\partial h-\bar\varphi=0$ and,
by Weyl's lemma, $g=h-\bar\psi$ is holomorphic on $\C.$
Since $g(z)\to0$ as $z\to\infty,$ we conclude that $g=0$ by the Liouville
theorem.
Thus, we have $h=\bar\psi,$ which is equivalent to the required relation.
\end{pf}

Let $c\ge0$ and $0\le r<1.$
Wirths \cite[Satz 1]{Wirths78} gave the sharp upper bound of $|g'(z)|,$
$|z|=r,$ in terms of $c$ and $r$ for holomorphic functions $g$ on $\D$ with
$\|g\|_c\le 1.$
His result specialized for $c=2$ and $r=0$ can be stated as follows.

\begin{lem}\label{lem:W}
Every holomorphic function $g$ on the unit disk $\D$ with $\|g\|_2\le 1$
satisfies
$$
|g'(0)|\le \frac{25\sqrt5}{16},
$$
where equality holds when $g(z)=(25\sqrt5/16)z.$
\end{lem}

\begin{pf}[Proof of Proposition \ref{prop:main}]
As in the proof of Theorem \ref{thm:main}, we can reduce the proof to the
assertions $|Vf(0)|\le (25\sqrt5/16) \|Sf\|_2$ and
$|Sf(0)|\le (4/3)\|Vf\|_3.$

Since $Vf(0)=(Sf)'(0),$
the first inequality and its sharpness follow from Lemma \ref{lem:W}.
Indeed, a function $f$ with $Sf(z)=cz$ for a constant $c$ 
(e.g.~$f(z)=e^{az}$ for a constant $a\ne0$) satisfies the equality.

We next show the second inequality.
By Lemma \ref{lem:int}, we have
\begin{align*}
|Sf(0)|
&\le\frac1{\pi}\iint_{|\zeta|<1}
\frac{(1-|\zeta|^2)^4|Vf(\zeta)|}{|\zeta|}d\xi d\eta \\
&\le\frac1{\pi}\iint_{|\zeta|<1}
\frac{(1-|\zeta|^2)\|Vf\|_3}{|\zeta|}d\xi d\eta \\
&=2\|Vf\|_3\int_0^1(1-r^2)dr=\frac43\|Vf\|_3.
\end{align*}
The proof is now complete.
\end{pf}

\begin{example}\label{ex:S}
Meromorphic functions $f$ on the unit disk with the property
that $Sf(z)=c(1-z^2)^{-2}$ for a constant $c$ are sometimes very important
(see for example \cite{CO93}).
The function $l(z)=\log\frac{1+z}{1-z}$ has the least norm $\|Sl\|_2=2$
within those which have no quasiconformal extension to the Riemann sphere,
and it satisfies $Sl(z)=2(1-z^2)^{-2}.$
Indeed, it is known that $f$ is never univalent in $\D$ if
$Sf(z)=c(1-z^2)^{-2}$ for a constant $c>2.$
On the other hand, the Koebe function $ k(z)=z/(1-z)^2$ has the
maximal norm $\|Sk\|_2=6$ within univalent meromorphic functions on $\D$
and again satisfies $Sk(z)=-6(1-z^2)^{-2}.$
By the following lemma, we see that $\|Vf\|_3=8\sqrt3 |c|/9$
whenever $Sf(z)=c(1-z^2)^{-2}$ for a constant $c.$
In particular, $\|Vl\|_3=16\sqrt3/9$ and $\|Vk\|_3=16\sqrt3/3.$
\end{example}

\begin{lem}
Suppose that $Sf(z)=(1-z^2)^{-2}$ in $z\in\D.$
Then $\|Vf\|_3=8\sqrt3/9.$
\end{lem}

\begin{pf}
By \eqref{eq:Vf}, we have
$$
Vf(z)=\frac{4z}{(1-z^2)^3}-\frac{4\bar z}{1-|z|^2}\cdot\frac1{(1-z^2)^2}
=\frac{8i\Im z}{(1-|z|^2)(1-z^2)^3}.
$$
Fix $r\in(0,1)$ for a while and set $z=re^{i\theta}.$
Then
$$
|Vf(z)|=\frac8{1-r^2}
\left(\frac{r^2\sin^2\theta}{(1-2r^2\cos2\theta+r^4)^3}\right)^{1/2}
=\frac8{1-r^2}F(r,\theta)^{1/2}.
$$
Since $F(r,\theta)=F(r,-\theta)=F(r,\pi-\theta),$
it is enough to consider the case when $0\le\theta\le\pi/2.$
We now have
$$
\frac{\partial F}{\partial\theta}
=\frac{2r^2\cos\theta\sin\theta(1-6r^2+r^4+4r^2\cos2\theta)}%
{(1-2r^2\cos2\theta+r^4)^4}.
$$
Note that $1-6r^2+r^4+4r^2\cos2\theta=0$ if and only if
$\sin\theta=(1-r^2)/(2\sqrt2 r)$
and that there exists $\theta=\theta(r)\in(0,\pi/2]$ satisfying this relation
only when $\sqrt3-\sqrt2\le r<1.$
In this case, $F(r,\theta)$ takes its maximum value at $\theta=\theta(r)$
within $0\le\theta\le\pi/2.$
Since $1-2r^2\cos(2\theta(r))+r^4=3(1-r^2)^2/2,$ we have
$$
F(r,\theta(r))=\frac{(1-r^2)^2/8}{27(1-r^2)^6/8}=\frac1{27(1-r^2)^{4}}.
$$
Thus
$$
(1-r^2)^3\max_{|z|=r}|Vf(z)|=\frac{8(1-r^2)^3}{3\sqrt3(1-r^2)^3}
=\frac8{3\sqrt3}
$$
for $\sqrt3-\sqrt2\le r<1.$
When $0<r<\sqrt3-\sqrt2,$ $F(r,\theta)$ is increasing in $0<\theta<\pi/2$
and thus $F(r,\theta)\le F(r,\pi/2)=r^2/(1+r^2)^6<1/[27(1-r^2)^{4}].$
Therefore, we conclude that $\|Vf\|_3=8/(3\sqrt3)=8\sqrt3/9.$
\end{pf}

\def\cprime{$'$} \def\cprime{$'$} \def\cprime{$'$}
\providecommand{\bysame}{\leavevmode\hbox to3em{\hrulefill}\thinspace}
\providecommand{\MR}{\relax\ifhmode\unskip\space\fi MR }
\providecommand{\MRhref}[2]{%
  \href{http://www.ams.org/mathscinet-getitem?mr=#1}{#2}
}
\providecommand{\href}[2]{#2}

\end{document}